\documentclass [11pt]{article}
\usepackage{amsfonts}
\usepackage{amsmath}
\usepackage{latexsym}
\usepackage{epsfig}
\usepackage{graphicx}
\def\e{\varepsilon}

\def\d{\delta}

\def\l{\lambda}

\newtheorem{theorem}{Theorem}[section]
\newtheorem{lemma}[theorem]{Lemma}
\newtheorem{corollary}[theorem]{Corollary}

\newtheorem{propos}[theorem]{Proposition}

\newtheorem{defin}[theorem]{Definition}

\newcommand{\proofstart}{{\bf Proof.\hspace{1.5em}}}

\newcommand{\proofend}{\hspace*{\fill}\mbox{$\Box$}}

\date{}

\begin{document}
\title{ Embedding nearly-spanning bounded degree trees}

\author{
Noga Alon \thanks{
Schools of Mathematics and Computer Science, Raymond and Beverly
Sackler Faculty of Exact Sciences, Tel Aviv University, Tel Aviv
69978, Israel and IAS, Princeton, NJ 08540, USA. Email:
nogaa@post.tau.ac.il. Research supported in part by a USA-Israeli BSF 
grant, by
NSF grant CCR-0324906, by a Wolfensohn fund and by the State of New
Jersey.}
\and
Michael Krivelevich \thanks{Department of Mathematics, 
Raymond and Beverly Sackler Faculty of Exact
Sciences, Tel Aviv University, Tel Aviv 69978, Israel. E-mail: 
krivelev@post.tau.ac.il. Research supported in part by USA-Israel
BSF Grant 2002-133, and by grants 64/01 and 526/05 from the Israel Science
Foundation.}
\and  Benny Sudakov \thanks{Department of Mathematics, Princeton 
University,
Princeton, NJ 08544. E-mail: bsudakov@math.princeton.edu.
Research supported in part by NSF CAREER award DMS-0546523, NSF grant
DMS-0355497, USA-Israeli BSF grant, and by an Alfred P. Sloan fellowship.} }
\maketitle
\begin{abstract}
We derive a sufficient condition for a sparse graph $G$ on
$n$  vertices to contain a copy of a tree $T$ of maximum
degree at most $d$ on $(1-\epsilon)n$ vertices, in
terms of the expansion properties of $G$. As
a result we show that for fixed $d\ge 2$ and 
$0<\epsilon<1$, there exists a constant $c=c(d,\epsilon)$
such that a random graph $G(n,c/n)$ contains almost surely
a copy of every tree $T$ on $(1-\epsilon)n$ vertices with maximum degree 
at most $d$. We also prove that if
an $(n,D,\lambda)$-graph $G$ (i.e., a $D$-regular graph
on $n$ vertices all of whose eigenvalues, except the first 
one, are at most $\lambda$ in their absolute values) has 
large enough spectral gap $D/\lambda$ as a function of $d$
and $\epsilon$, then $G$ has a copy of every tree $T$ as
above. 
\end{abstract}

\section{Introduction}
In this paper we obtain a sufficient condition for a sparse graph $G$ 
to contain a copy of every nearly-spanning tree $T$ of bounded 
maximum degree, in terms of the expansion properties of $G$.
The restriction on the degree of $T$ comes naturally from the fact that we consider graphs 
of constant degree. Two important examples where our condition 
applies are random graphs
and graphs with a large spectral gap. 

The random graph $G(n,p)$ denotes the probability space whose
points are graphs on a fixed set of $n$ vertices, where each pair
of vertices forms an edge, randomly and independently, with
probability $p$. We say that the random graph $G(n,p)$
possesses a graph property $\cal P$ {\em almost surely},
or a.s. for short, if the probability that $G(n,p)$ satisfies $\cal P$
tends to 1 as the number of vertices $n$ tends to infinity.

The problem of existence of large trees with specified shape in random graphs 
has a long history with most of the results being devoted to finding a long path.
Erd\H{o}s conjectured that a random graph $G(n,c/n)$ a.s. contains a path of length
at least $(1-\alpha(c))n$, where $\alpha(c)$ is a constant smaller than one for
all $c>1$  and $\lim_{c \rightarrow \infty} \alpha(c)=0$. This conjecture was proved 
by Ajtai, Koml\'os and Szemer\'edi \cite{AKS} and, in a slightly weaker form, by 
Fernandez de la Vega \cite{FD1}. These results were significantly improved by 
Bollob\'as \cite{B}, who showed that $\alpha(c)$ decreases exponentially in $c$.
Finally Frieze \cite{F} determined the correct speed of convergence of 
$\alpha(c)$ to zero and proved that $\alpha(c)=(1+o(1))ce^{-c}$. 

The question of existence of large trees of bounded degree other than paths in sparse random graphs was studied
by Fernandez de la Vega in \cite{FD2}. He proved that there exist two constants $a_1>0$ and $a_2>0$ such that for fixed tree
$T$ of order $n/a_1$ with maximum degree at most $d$ the random graph $G(n, c/n)$ 
with $c=a_2d$ almost surely contains $T$. The constant $a_1$ in this result is rather large and allows to 
embed only trees that occupy a small proportion of the random 
graph. Also, observe that Fernandez de la Vega's result gives the almost
sure existence of a {\em fixed} tree $T$, and not of all such
trees simultaneously.
  
Our first theorem improves the result of 
Fernandez de la Vega and generalizes the above mentioned results 
on the existence of long 
paths. It shows that the sparse random graph contains almost 
surely every nearly-spanning tree of bounded degree. 

\begin{theorem}
\label{random}
Let $d\geq 2$, $0<\e<1/2$ and let
$$c \geq \frac{10^6 d^3 \log d \log^2 (2/\e)}{\e}.$$
Then almost surely the random graph $G(n, c/n)$ contains every tree 
of maximum degree at most $d$
on $(1-\e)n$ vertices.
\end{theorem}

Results guaranteeing the existence of a long path in a sparse
graph can be obtained in a more general situation when the host
graph has certain expansion properties. Given a graph
$G=(V,E)$ and a subset $X \subset V$ let
$N_G(X)$ denote the set of all neighbors of vertices of $X$ in 
$G$. Using his celebrated 
rotation-extension technique, P\'osa \cite{P} proved that if for every $X$ in $G$ with $|X|\leq k$
we have that $|N_G(X)\setminus X| \geq 2|X|-1$, then $G$ contains a path of length $3k-2$.
A remarkable generalization of this result from paths to trees of bounded degree was obtained by 
Friedman and Pippenger \cite{FP}. They proved that if 
$|N_G(X)| \geq (d+1)|X|$ for every set $X$ in $G$ with $|X|\leq 2k-2$, then 
$G$ contains every tree with $k$ vertices and maximum degree at most $d$.
Note that this result allows to embed only trees whose size is relatively small 
compared to the size of $G$. What if we want to embed trees which are 
nearly-spanning?
It turns out that a slightly stronger expansion property, based on the spectral gap condition, 
is already enough to attain this goal. 

For a graph $G$ let $\lambda_1 \geq \lambda_2 \geq \ldots \geq \lambda_n$
be the eigenvalues of its adjacency matrix. The quantity $\lambda(G)=
\max_{i \geq 2} |\lambda_i|$ is called the {\em second eigenvalue} of $G$.
A graph $G=(V,E)$ is called an {\em $(n,D,\lambda)$-graph} if it is
$D$-regular, has $n$ vertices and the second eigenvalue of $G$ is at most
$\lambda$. It is well known (see, e.g., \cite{AloSpe00} for more details) that if $\lambda$ is
much smaller than the degree $D$, then $G$ has strong expansion properties, so the ratio
$D/\lambda$ could serve as some kind of measure of expansion of $G$.
Our next result shows that an $(n,D,\lambda)$-graph $G$ with large enough spectral gap $D/\lambda$ 
contains a copy of every nearly-spanning tree with bounded degree.

\begin{theorem}
\label{pseudo-random}
Let $d\geq 2$, $0<\e<1/2$ and let $G$ be an $(n,D,\l)$-graph such
that
$$\frac{D}{\lambda} \geq  \frac{160d^{5/2}\log (2/\e)}{\e}.$$
Then $G$ contains a copy of every tree $T$ with $(1-\e)n$ vertices
and with maximum degree at most $d$.
\end{theorem}

Our main results are tightly connected to the notion of universal
graphs. For a family ${\cal H}$ of graphs, a graph $G$ is ${\cal
H}$-universal if $G$ contains every member of ${\cal H}$ as a
(not necessarily induced) subgraph. The construction of sparse
universal graphs for various families arises in several fields
such as VLSI circuit design, data representation and parallel
computing (see, e.g., the introduction of \cite{ACKRRS} for a 
short survey and relevant references). Our two main results show
that sparse random graphs and pseudo-random graphs on $n$
vertices are universal
graphs for the family of bounded-degree trees on
$(1-\e)n$ vertices.
Quite an extensive research exists on universal graphs for trees
\cite{BCLR}, \cite{CG1}, \cite{CG2}, \cite{CG3}, \cite{CGP},
\cite{FP}. 
The most interesting result is that of Bhatt et al. who
showed in \cite{BCLR} that there exists a universal graph $G$ on
$n$ vertices for the family of trees on $n$ vertices with maximum
degree $d$, whose maximum degree is bounded by a function of $d$.
It is
instructive to compare our results with those of \cite{BCLR}:
they succeed in embedding spanning trees as opposed to 
nearly-spanning in our case; on the other hand, their universal graph is
a concrete carefully constructed graph that has very dense pieces
locally, while we provide a very large family of universal graphs
possessing many additional properties that can be useful 
for obtaining further results on universal graphs.  

The results of Theorem \ref{random} and \ref{pseudo-random} 
can be deduced from a more general statement
which we present next. We need the following definition.

\begin{defin}\label{def1}
Given two positive numbers $c$ and $\alpha<1$, a graph $G=(V,E)$ is 
called an $(\alpha,c)$-expander
if every subset of vertices $X\subset V(G)$ of size 
$|X|\le \alpha |V(G)|$ satisfies:
$$
|N_G(X)|\ge c |X|\ .
$$
\end{defin}

\begin{theorem}\label{th1}
Let $d\ge 2$, $0<\e<1/2$. 
Let $G=(V,E)$ a graph on $n$ vertices of minimum
degree $\delta=\delta(G)$ and maximum degree $\Delta=\Delta(G)$. Let 
$n, \delta, \Delta$ satisfy:
\begin{enumerate}
\item (order of graph is large enough)
$$n \ge \frac{480d^3\log (2/\e)}{\e};$$
\item (maximum degree is not too large compared to the minimum degree)
$$\Delta^2 \le \frac{1}{K} e^{\delta/(8K)-1}~~ 
\mbox{where} ~~K=\frac{20 d^2 \log (2/\e)}{\e}.$$
\item (local expansion)
Every induced subgraph $G_0$ of $G$  
with minimum degree at least $\frac{\e \delta}{40d^2 \log (2/\e)}$
is a $(\frac{1}{2d+2},d+1)$-expander.
\end{enumerate}
Then $G$ contains a copy of every tree $T$ on at most $(1-\e)n$ vertices
of maximum degree at most $d$.
\end{theorem}  

The rest of this paper is organized as follows. In the next two sections we 
show how Theorem \ref{th1} can be used to embed nearly-spanning 
trees of bounded 
degree in random and pseudo-random graphs. We present the proof of 
Theorem \ref{pseudo-random} first, since it is short and less technical, and then prove Theorem 
\ref{random}. In Section 4 we describe the plan of the proof 
of Theorem \ref{th1} and discuss some 
technical tools needed to fulfill this plan. 
The proof of this theorem appears in Section 5. The last section of the paper contains several 
concluding remarks and open problems. 

Throughout the paper we make no attempts to optimize the
absolute constants. To simplify the presentation, we often omit
floor and ceiling signs whenever these are not crucial.
Throughout the paper, $\log$ denotes logarithm in the natural base $e$.

\section{Embedding in pseudo-random graphs}
In this section we prove Theorem 
\ref{pseudo-random}. First we need the following lemma that
shows that an $(n,D,\l)$-graph has 
the local expansion property required by 
Condition 3 of Theorem \ref{th1}.

\begin{lemma}\label{prg}
Let $d\ge 2$. Let $G=(V,E)$ be an $(n,D,\l)$-graph. Denote 
$$
D_0=\frac{2\l (d+1)}{\sqrt{d}}\ .
$$
Then every induced subgraph $G_0$ of $G$ of 
minimum degree at least $D_0$ is a
$\big(\frac{1}{2d+2},d+1\big)$-expander.
\end{lemma}

\proofstart
We will use the following well known estimate on the edge distribution
of an $(n,D,\l)$-graph $G$ (see, e.g., \cite{AloSpe00},
Corollary 9.2.5).
For every two (not necessarily disjoint)
subsets $B,C\subseteq V$, let $e(B,C)$ denote  
the number of ordered pairs $(u,v)$ with $u \in B, v \in C$ such that
$uv$ is an edge. Note that if $u,v \in B \cap C$, then the edge
$uv$ contributes $2$ to $e(B,C)$. In this notation, 
$$
\left|e(B,C)-\frac{|B||C|D}{n}\right|\le \l\sqrt{|B||C|}\ .
$$

Let $U$ be a subset of vertices of $G$ such that the induced subgraph
$G_0=G[U]$ has minimum degree at least $D_0$. 
Suppose that the claim is false. Then there exists a subset 
$X\subset U$ of size $|X|=t \leq |U|/(2d+2)$ satisfying 
$|N_{G_0}(X)|< (d+1)|X|$. 
  
By the above estimate with $B=X$ and $C=N_{G_0}(X)$  we have:
$$
D_0t\le e(B,C)\le \frac{t(d+1)tD}{n}
+\l t \sqrt{ d+1}\,
$$
and therefore
\begin{equation}\label{eq1}
\frac{t}{n} \ge \frac{D_0}{(d+1)D}-\frac{\l}{\sqrt{d+1}D}\ .
\end{equation}
Also, note that there are no edges of $G$ from $X$ to
$Y=U-(X \cup N_{G_0}(X))$ as $G_0$ is an induced subgraph of $G$.
From $t=|X|\le |U|/(2d+2)$ and $|N_{G_0}(X)| \leq (d+1)t$ it follows 
that $|Y|\ge dt$. Thus
$$
0=e(X,Y)\ge \frac{t (dt) D}{n}-\l \sqrt{t (dt)}\,,
$$
implying 
\begin{equation}\label{eq2}
\frac{t}{n}\le \frac{\l}{\sqrt{d} D}\ .
\end{equation}
Comparing (\ref{eq1}) and (\ref{eq2}) we obtain
$$
\frac{D_0}{(d+1)D}-\frac{\l}{\sqrt{d+1}D}\le \frac{\l}{\sqrt{d}D}\ .
$$
Plugging in the definition of $D_0$ we derive a contradiction.
\proofend

\vspace{0.25cm}
\noindent
{\bf Proof of Theorem \ref{pseudo-random}.}\hspace{1.5em}
Since every graph with minimal degree $k$ contains all the trees on $k$ vertices
we can assume that $D \leq (1-\e)n$. Let $A$ be the adjacency matrix of $G$. The trace of $A^2$ 
equals the number of ones in $A$, which is exactly
$2|E(G)|=nD$. We thus obtain that
$$
nD=Tr(A^2)=\sum_{i=1}^n\lambda_i^2\le D^2+(n-1)\lambda^2
$$
and therefore $\lambda^2 \geq \frac{D(n-D)}{n-1} \geq \e D$.
This together with our assumption on $D/\lambda$ implies
$$n \geq D \geq \e\left(\frac{ D}{\lambda}\right)^2 \geq 
\frac{160^2 d^5 \log^2 (2/\e)}{\e}.
$$
Since $\Delta(G)=\d(G)=D$, from this inequality 
it follows that
$G$ satisfies Conditions 1, 2 of Theorem \ref{th1}. 
Finally since
$$ \frac{\e D}{40 d^2 \log (2/\e)}
\geq \frac{2(d+1)\l}{\sqrt{d}}$$
we can conclude using Lemma \ref{prg} that $G$ also 
satisfies the last condition of Theorem \ref{th1}.
Thus $G$ contains every tree of size $(1-\e)n$ with
maximum degree at most $d$.
\hfill $\Box$

\section{Embedding in random graphs}
To prove Theorem \ref{random} we first need to show that a sparse random graph contains a.s. a nearly
spanning subgraph with good local expansion properties.

\begin{lemma}
\label{rg}
For every integer $d\ge 2$, real $0<\theta<1/2$ and $D \geq 50\theta^{-1}$
the random graph $G\big(n,\frac{4D}{n}\big)$ almost surely contains a subgraph $G^*$ having the
following properties:
\begin{enumerate}
\item $|V(G^*)|\ge (1-\theta)n$;
\item $D\le d_{G^*}(v)\le 10D$ for every $v\in V(G^*)$;
\item every induced subgraph $G_0$ of $G^*$ of minimum degree at least
$D_0=100d\log D$ is a $\big(\frac{1}{2d+2},d+1\big)$-expander. 
\end{enumerate}
\end{lemma}

The following statement contains a few easy facts about random graphs.
\begin{propos}
\label{properties}
Let $G(n,p)$ be a random graph with $np>20$, then almost surely

\noindent
$(i)$\, The number of edges between any
two disjoint subsets of vertices $A, |A|=a$ and $B, |B|=b$
with $abp \geq 32n$ is at least $abp/2$ and at most $3abp/2$.

\noindent
$(ii)$\, Every subset of vertices of size $a \leq n/4$ spans
less than $anp/2$ edges.
\end{propos}

\proofstart
(i)\, Since the number of edges between $A$ and $B$ is a
binomially distributed random variable with parameters $ab$ and $p$, it follows by the
standard Chernoff-type estimates (see, e.g., \cite{AloSpe00}) that
(denoting $t=abp/2$)
$$\mathbb{P}\Big[e(A,B)-abp < -t\Big] \leq
e^{-\frac{t^2}{2abp}}=e^{-abp/8}$$
and
$$\mathbb{P}\Big[e(A,B)-abp > t\Big] \leq
e^{-\frac{t^2}{2abp}+\frac{t^3}{2(abp)^2}}=e^{-abp/16}.$$
Using that $abp\geq 32n$ we can bound the probability that there are sets $A, B$ with
$|e(A,B)-abp|>abp/2$ by $2^n \cdot 2^n \cdot \big(2e^{-2n}\big)=o(1)$.

\noindent
Since $np/2 \geq 10$ and $n/a \geq 4$, the probability
that there is a subset of size $a$ which 
violates the assertion (ii) is at most
$$\mathbb{P}_a \leq {n \choose a} {a^2/2 \choose anp/2}p^{anp/2} \leq 
\left(\frac{en}{a}  \Big(\frac{ea}{np}\Big)^{np/2} p^{np/2}\right)^a=
\left(\frac{e^{np/2+1}}{(n/a)^{np/2-1}}\right)^a \leq 
\left(\frac{e^{11}}{(n/a)^9}\right)^a.$$ 
It is easy to see that $\mathbb{P}_a \ll n^{-1}$ for all $a \leq n/4$ and so
$\sum_a\mathbb{P}_a=o(1)$.
\proofend

\vspace{0.25cm}
\noindent
{\bf Proof of Lemma \ref{rg}.}\hspace{1.5em}
Let $G=G(n,p)$ be a random graph with $p=\frac{4D}{n}$ and let $X$
be the set of $\theta n/2$ vertices of largest degrees in $G$.
By Part (ii) of 
Proposition  \ref{properties}, a.s. this set spans less than 
$|X|np/2=2D|X|$ edges. Also, since
$\frac{4D}{n}|X|(n-|X|) \geq 2D \theta (n/2)\geq 50 n$, Part (i) of this proposition implies 
that a.s. the number of edges between $X$ and $V(G)-X$ is at most 
$3|X|np/2=6D|X|$. Therefore the sum of the 
degrees of the vertices in $X$ is bounded by $10D|X|$
and hence there is a vertex in $X$ with degree at most 
$10D$. By definition of $X$, this implies that there are 
at most $\theta n/2$ vertices in $G$ with degree larger than $10D$. Delete 
these vertices and denote the remaining graph by $G'$. Next as long as 
$G'$ contains a vertex $v$ 
of degree less than $D$, delete it. If we deleted more than $\theta n/2$ vertices, then the
original random graph contains two sets $Y$ and $V(G')-Y$ such that 
$|Y|=\theta n/2$, $|V(G')-Y|\geq (1-\theta)n\geq n/2$ and there are less than
$D|Y|\leq p|Y||V(G')-Y|/2$ edges between them. Since $\frac{4D}{n}|Y||V(G')-Y| 
\geq \theta D n\geq 50 n$, again by Part (i) of the previous 
statement this a.s. does not happen. Denote the resulting graph by $G^*$.
Then it satisfies the first two conditions of the lemma and it remains to verify the third 
condition.

Suppose to the contrary that $G^*$ contains a 
subset of vertices $U$ such that the 
induced subgraph $G_0=G^*[U]$ has minimum degree 
at least $D_0=100d\log D$
and is not a $\big(\frac{1}{2d+2},d+1\big)$-expander. 
Then there exists a set $X\subset U$ of size $|X|=t$ 
such that the set $C= N_{G_0}(X)$ has size at most
$(d+1)t$ and there are at least $D_0|X|/2=50dt\log D$ edges
with an end in $X$ and another end in $C$. 
If $t\le \frac{\log D}{D}n$, then the probability that $G(n,p)$
contains such sets is at most

\begin{eqnarray*}
\mathbb{P}_t &\leq& {n \choose t}{n\choose 
{(d+1)t}}{{{t(d+1)t}}\choose{50d t\log D}}
p^{50dt \log D }\\
&\le&
\left[\left(\frac{en}{t}
\right) \left( \frac{en}{(d+1)t}\right)^{d+1}
      \left(\frac{e(d+1)tp}{50d\log D}\right)^{50d\log D}\right]^t\\
&\le&\left[e\,\left(\frac{n}{t}\right)^{2d}
      \left(\frac{e}{8}\cdot\frac{Dt}{n\log D}\right)^{50d \log D}
\right]^t\\
&=& \left[e\,\left(\frac{e}{8}\right)^{50d \log D}
 \left(\frac{D}{\log D}\right)^{2d} 
\left(\frac{Dt}{n\log D} \right)^{50d \log D-2d}\right]^t \\
&<&
\left[e^{-25d \log D+2d\log D+1}
\left(\frac{t}{n\log D/D} \right)^{40d \log D}\right]^t \\
&\leq&  
\left[D^{-20d}
\left(\frac{t}{n\log D/D} \right)^{40d\log D}\right]^t
\ .
\end{eqnarray*}  
Checking separately two cases $t < {\log n}$ and  
${\log n} \leq t < \frac{\log D}{D}n$ it is easy
to see that in both $\mathbb{P}_t \ll n^{-1}.$

If  $t\ge \frac{\log D}{D}n$ we apply a different argument. 
Note that there are no edges of $G(n,p)$ from $X$ to
$C=U-(X \cup N_{G_0}(X))$ since $G_0$ is an induced subgraph.
From $t=|X|\le |U|/(2d+2)$ and $|N_{G_0}(X)| \leq (d+1)t$ it follows
that $|C|\ge dt$ and therefore the probability of such event in $G(n,p)$
is at most
\begin{eqnarray*}
\mathbb{P}_t&\leq& {n\choose t}{{n}\choose{dt}}(1-p)^{dt^2}\le
\left[\frac{en}{t}\,\cdot\, \left(\frac{en}{dt}\right)^de^{-pdt}
\right]^t \\
&\leq& \left[\left(\frac{en}{t}\right)^{2d}e^{-pdt}\right]^t
=
\left[\left(\frac{en}{t}\right)^2\,\cdot\,e^{-pt}\right]^{dt}\\
&\leq& \left[\left(\frac{en}{n\log D/D}\right)^2\,\cdot\,
e^{-\frac{4D}{n}\,\cdot\,\frac{\log D}{D}n}\right]^{dt}\\
&\leq& \left( D^2 D^{-4}\right)^{dt}=o(n^{-1})\ . 
\end{eqnarray*}
Thus the probability that $G^*$ fails to satisfy the third condition is at most
$\sum_{t=1}^n \mathbb{P}_t=o(1)$.
\hfill $\Box$

\vspace{0.25cm}
\noindent
{\bf Proof of Theorem \ref{random}.} \hspace{1.5em}
Let $d \geq 2$, $0<\e<1/2$ and $c$ satisfy the assumption of Theorem \ref{random}.
Set $\theta=0.01\e$, $D=c/4$ and let 
$\e_1=\frac{\e-\theta}{1-\theta}\geq 0.99\e$.
Then by Lemma \ref{rg} $G(n, c/n)$ 
almost surely contains a subgraph $G^*$ of order
$n_1\geq (1-\theta)n$ such that $D \leq \d(G^*)\leq \Delta(G^*) \leq 10D$ and 
every induced subgraph of $G^*$ with minimum degree at least
$100d \log D$ is an $\big(\frac{1}{2d+2},d+1\big)$-expander.
Using that $\Delta(G^*) \leq 10\d(G^*)$ and 
$$n_1 \geq \d(G^*) \geq D \geq \frac{10^6d^3\log d \log^2 (2/\e)}{4\e}>
\frac{480 d^3 \log(2/\e_1)}{\e_1}$$
we conclude that  $G^*$ satisfies Conditions 1 and 2
of Theorem \ref{th1} (with $\e_1$). 
To verify the third condition it is enough to check that
the assumptions in Theorem \ref{random} imply that
$$100d\log D \leq \frac{\e_1 D}{40d^2 \log (2/\e_1)}.$$
Note that one can simply substitute the 
lower bound for $D$ in the above expression, since 
$x/\log x$ is an increasing function for $x>3$.
Therefore by Theorem \ref{th1}, $G^*$ contains every tree of size
$(1-\e_1)n_1 \geq (1-\e_1)(1-\theta)n=(1-\e)n$ with maximum degree at most $d$.
\hfill $\Box$

\section{Embedding plan and main tools}
To prove Theorem \ref{th1} we will use the following framework. 
Given a tree $T$, we first cut it into subtrees 
$T_1, T_2, \ldots ,T_s$ of 
carefully chosen sizes, so that the number $s$
of these subtrees 
satisfies $s \leq 10d^2 \log(2/\e)$, and each subtree $T_i$ is connected
by a unique edge to the union of all previous subtrees. The subtrees
$T_i$ will be embedded sequentially in order, starting from $T_1$.

We then choose $s$ pairwise disjoint sets of vertices  $S_1,S_2,
\ldots, S_s$ whose total size is at most $\e n/2$, such that each
vertex of the graph has many neighbors in each set $S_i$. The set
$S_i$ will be used only when embedding the subtree $T_i$, and 
will not be touched before that step.  (It can be used later,
but we will not do it here, as it complicates matters and does not
improve the estimates in any essential way). During the  embedding
process we maintain a set $R$ of at most $s$ vertices, which
will consist of all roots of the  trees $T_i$ that still
have to be embedded, and will not contain any vertex of the sets
$S_i$.

At the $i$-th step we are to embed the tree $T_i$ starting from a given
root $x_i\in R$. (At the first step a root is chosen arbitrarily.)
Suppose that the current set of unused vertices of $G$ is $V_{i-1}$.
We take an arbitrary subset $U_i$ of size $|U_i|=\Theta(|V(T_i)|d)$ 
which contains the  vertex $x_i$ that will be the
root  of $T_i$ (but contains no other members of $R$), 
contains the set $S_i$, and contains no member of $S_j$ for $j \neq i$.
Note that as each vertex has many neighbors in $S_i$, the 
minimum degree in the induced subgraph of $G$ on $U_i$ is large.

Since the minimum degree of $G[U_i]$ is large enough, we can
use the result of Friedman and Pippenger to embed a copy of $T_i$ 
in $U_i$, rooting
it at $x_i$. (We actually need a slightly modified, rooted version, of their result).
All vertices 
of $U_i$ unused when embedding $T_i$ are
recycled, and we thus get $V_i$ by deleting from $V_{i-1}$ only the
vertices used for embedding $T_i$.

The final step of embedding $T_i$ is to embed the edges crossing from
$T_i$ to yet unembedded pieces $T_j$, with $j>i$ (using vertices of $U_i$). 
We then add the endpoints of those
edges outside $T_i$ to the list $R$ of special vertices, and delete
$x_i$ from $R$. Each of the newly added special vertices will 
serve as a
root for embedding the corresponding piece $T_j$. Observe that the
number of special vertices is at most $s$ at any stage of 
the embedding.

The precise technical details are described in
what follows.

\subsection{The result of Friedman and Pippenger}
The cornerstone of our proof is the embedding result of Friedman
and Pippenger. In fact, we need a slightly stronger version of it
-- they showed the existence of a tree $T$, while we need to
embed a rooted version of $T$ in $G$ starting from a fixed vertex
$v\in V(G)$ as its root. Luckily, a careful reading of \cite{FP}
reveals that the following holds as well.

\begin{theorem}[\cite{FP}]\label{FP}
Let $T$ be a tree on $k$ vertices of maximum degree at most $d$ rooted
at $r$. Let $H=(V,E)$ be a non-empty graph such that, for every
$X\subset V(H)$ with $|X|\le 2k-2$,
$$
|N_H(X)|\ge (d+1) |X|\ .
$$
Let further $v\in V(H)$ be an arbitrary vertex of $H$. Then $H$ contains
a copy of $T$, rooted at $v$.
\end{theorem}

\subsection{Cutting the tree into pieces}
\begin{propos}\label{tree-cut1}
Let $d \geq 2$ and $k$ be positive integers. Let $T$ be a tree on at least $k+1$
vertices with maximum degree at most $d$. Then there exists an edge
$e\in E(T)$ such that at least one of the two trees obtained from $T$ by
deleting $e$ has at least $k$ and at most $(d-1)(k-1)+1$ vertices.
\end{propos}

\proofstart
Choose a leaf $r$ of $T$ arbitrarily and root $T$ at $r$. For $i\ge 0$
denote by $L_i$ the set of vertices of $T$ at distance $i$ from $r$.
For a vertex $v\in V(T)$ let $t(v)$ be the number of vertices in the
subtree of $T$ rooted at $v$. Now, let
$$
i_0=\max\{i: L_i\mbox{ contains a vertex $v$ with $t(v)\ge k$}\}\ .
$$
As $L_1$ has only one vertex $v$ with $t(v)=|V(T)|-1\ge k$, it follows
that $i_0\ge 1$. Choose a vertex $u\in L_{i_0}$ such that $t(u)\ge k$.
Then by the definition of $i_0$ all sons $w$ of $u$ in $T$ satisfy:
$t(w)\le k-1$, the number of sons does not exceed $d-1$, and therefore
$t(u)\le (d-1)(k-1)+1$. 

Let now $x$ be the father of $u$ in $T$. Then $e=(x,u)$ is the required edge.
\proofend 

\begin{corollary}\label{tree-cut2}
Suppose $0 < \e <1/2$, and let $T$ be an arbitrary tree 
on $(1-\e)n$ vertices, with maximum degree at most $d$. 
Then one can cut $T$ into $s$ subtrees
$T_1,T_2, \ldots ,T_s$, so that each tree $T_i$ is connected
by a unique edge to the union of all trees $T_j$ with $j<i$,
and such that for every $i>1$, 
$$
\frac{\e n/2 +\sum_{j>i} |V(T_j)|}{8d^2}
\leq   |V(T_i)| \leq 
\frac{\e n/2 +\sum_{j>i} |V(T_j)|}{8d}.
$$ 
For $i=1$,  the upper bound holds, but the lower bound may fail.
Moreover, $s \leq 10d^2 \log (2/\e).$
\end{corollary}

\proofstart 
We choose the trees $T_i$ one by one, starting from the 
{\em last} one.  By 
Proposition \ref{tree-cut1}  we first find a tree $T'_1$ of size
at least $\frac{\e n}{16d^2}$ and at most $\frac{\e n}{16d}$,
and omit it from $T$. Suppose we have already chosen
$T'_1, T'_2, \ldots ,T'_{i-1}$ such that for every $j<i$
$$
\frac{\e n/2 +\sum_{r<j} |V(T'_r)|}{8d^2} \leq  |V(T'_j)| \leq 
\frac{\e n/2 +\sum_{r<j} |V(T'_r)|}{8d},
$$ 
and each $T'_j$ has a unique edge joining it to 
$V(T)\setminus \cup_{r<j}V(T'_r)$.
Let $T'$ be the tree obtained from $T$ by omitting all the vertices
of all subtrees $T'_j$, $j<i$. If the number of vertices
of $T'$ is at most 
$\frac{\e n/2 +\sum_{r<i} |V(T'_r)|}{8d},$ then define 
$T'_i=T'$ and $s=i$. Else, apply Proposition \ref{tree-cut1}
to find a tree $T'_i$ in $T'$ whose size is at least
$\frac{\e n/2 +\sum_{r<i} |V(T'_r)|}{8d^2}$ and at most
$\frac{\e n/2 +\sum_{r<i} |V(T'_r)|}{8d}$, and continue.
To estimate the number of steps in this process define
$a_i=\e n/2+\sum_{j\leq i}|V(T'_j)|$. Observe that $a_0=\e n/2, 
a_i\leq \e n/2+|V(T)|\leq n$ and 
$a_i=a_{i-1}+|V(T'_i)| \geq \left(1+\frac{1}{8d^2}\right)a_{i-1}$. From this it follows that 
$$\frac{2}{\e} \geq \frac{a_i}{a_0} \geq \left(1+\frac{1}{8d^2}\right)^i,$$
and hence this process terminates after at most $10 d^2 \log (2/\e)$ steps.
Finally, for $1\le i\le s$ define $T_i=T'_{s-i+1}$.
\proofend 

\subsection{Splitting vertex degrees}
\begin{lemma}\label{l44}
Let numbers $K, \delta,\Delta$ satisfy
$$
K\Delta^2 e^{-(\delta/8K)+1} <1.
$$
Then the following holds.
Let $H=(V,E)$ be a graph in which $\delta \le d(v)\le \Delta$ for each
$v\in V$. 
Then $H$ contains $K$ pairwise disjoint sets of vertices
$S_1, S_2, \ldots ,S_K$ such that every vertex of $H$ has at least
$\frac{\delta }{2K}$ neighbors in each set $S_i$.
\end{lemma}

\proofstart 
This is a simple consequence of the Lov\'asz Local Lemma  (c.f., e.g., 
\cite{AloSpe00}, Chapter 5). Color the vertices of $H$ randomly and
independently by $K$ colors. For each vertex $v$ and color
$i$, $1 \leq i \leq K$, let $A_{v,i}$ be the event that $v$ 
has less than $\frac{\delta }{2K}$ neighbors  of color $i$.
By Chernoff's Inequality  the probability of each event
$A_{v,i}$ is at most $e^{-\delta/(8K)}$. In addition, each event
$A_{v,i}$ is mutually independent of all events but the events
$A_{u,j}$ where either $u=v$ or $u$ and $v$ have common neighbors in
$H$. As there are less than $K\big(\Delta(\Delta-1)+1\big)\leq K\Delta^2$ such events, 
it follows by the Local Lemma that with positive probability none
of the events $A_{v,i}$ holds. The desired result follows, by letting
$S_i$ denote the set of all vertices of color $i$.
\proofend

\section{Proof of Theorem \ref{th1}}
Let $G=(V,E)$ be a graph satisfying the assumptions of the theorem.
Let $T$ be a tree on at most $(1-\e)n$ vertices with maximum
degree at most $d$. By Corollary \ref{tree-cut2} the tree $T$ can be
partitioned into subtrees $T_1, T_2, \ldots, T_s$ satisfying the conditions 
of the Corollary, where $s \leq 10 d^2 \log (2/\e).$ 
Choose an arbitrary root for  $T_1$. For $i>1$, the root of $T_i$
is the vertex incident with the unique edge that connects $T_i$
to the union of the previous trees.
Put $K=2s/\e$.
By Condition 2 in the theorem, and Lemma \ref{l44} there are $K$ pairwise
disjoint sets of vertices $S_i$ of $G$ such that every vertex of $G$ has 
at least $ \delta/(2K) \geq \frac{\e \delta}{40 d^2 \log (2/\e)}$
neighbors in each set $S_i$. Take the $s$ smallest sets $S_i$ and renumber them
so that they are denoted by $S_1,S_2,\ldots ,S_s$. Obviously, their total
size is at most $\frac{ns}{K} =\frac{\e n}{2}$. We will not use all the other sets
$S_i, i>s$ in the rest of the proof.

Let $x_1$ be an arbitrary vertex of $G$ that does not lie in any of the sets
$S_i$. Define $R=\{x_1\}$ and let $U_1$ denote the set of all vertices
of $G$ besides those in $\cup_{j \neq 1} S_j$. As $U_1$ contains $S_1$, every
vertex in the induced subgraph $G[U_1]$ of $G$ on $U_1$
has degree at least  $\frac{\e \delta}{40 d^2 \log (2/\e)}$. Therefore,
by Condition 3 in Theorem \ref{th1}, and by Theorem \ref{FP}
there is a copy of $T_1$ in $G[U_1]$ rooted at $x_1$. (Note that 
by  Corollary \ref{tree-cut2}, 
the size of $T_1$ is at most $\frac{|U_1|}{8d}$ and hence indeed
one can apply here Theorem \ref{FP}.) Moreover, we can in fact embed
in $U_1$ the required tree $T_1$ together with the edges connecting it to 
the trees $T_j$ with $j>1$. Add the endpoints of these edges
to the list $R$ of planned roots for the trees $T_j$, and delete
$x_1$ from $R$. This completes the embedding of $T_1$.

Assume that we have already embedded the first $i-1$ trees $T_r$,
where each $T_r$ has been rooted in the vertex of $R$ specified as its
root, and where in step number $r$ the tree $T_r$ has been embedded
using no vertices of $R$ besides its root, and no vertices of
$\cup_{j \neq r} S_j$, we proceed to the $i$-th step, in which
we are to embed the tree $T_i$ starting from a given root
$x_i \in R$. Let $U_i$ be the set of all vertices  of $G$ that have not
been used for embedding the part of $T$ embedded so far, 
besides the vertices in $R-\{x_i\}$ and besides the
vertices in $\cup_{j \neq i} S_j$. 
As before, since
$U_i$ contains $S_i$, every
vertex in the induced subgraph $G[U_i]$ of $G$ on $U_i$
has degree at least  $\frac{\e \delta}{40 d^2 \log (2/\e)}$. Therefore,
by Condition 3 in Theorem \ref{th1}, it is a
$(\frac{1}{2d+2},d+1)$-expander. By Corollary \ref{tree-cut2}, the size of $T_i$ is at most 
$|U_i|/8d$ and the number of edges connecting $T_i$ to the trees 
$T_j$ for $j>i$ is bounded by
$$s \leq 10 d^2 \log (2/\e) \leq \frac{\e n}{48d} \leq \frac{|U_i|}{24d}.$$
Therefore the
size of $T_i$ together with the vertices in $T_j$ for $j>i$ that
are connected to it is less than a fraction $\frac{1}{6d} \leq \frac{1}{2(2d+2)}$ of
the size of $U_i$. 
Therefore, by Theorem \ref{FP} we can embed $T_i$ including the edges
connecting it to the trees $T_j$ with $j>i$ in $U_i$,
rooting it at $x_i$, and add the endpoints of the edges from $T_i$
to future $T_j$'s to $R$. 

As this process can be  carried out until we finish the embedding of
$T_s$, the assertion of the theorem follows.
\proofend

\section{Concluding remarks}
\begin{itemize}
\item
Our lower bound on the edge probability of a random 
graph in Theorem \ref{random}
seems far from being best possible, and the correct order of magnitude 
should probably be more similar 
to the case of a longest path. Hence it is likely that already when 
$c=O(d\log (1/\epsilon))$ the random graph $G(n,c/n)$ contains a.s.
every tree on 
$(1-\epsilon)n$ vertices with maximum  degree at most $d$.
\item
Embedding {\em spanning} trees of bounded degree in sparse random graphs is an 
intriguing question which is completely open. In case of the 
path, this question is very well understood
(see, e.g., Chapter 8 of \cite{Bol}) and it is known that for 
$p=O(\log n/n)$ the random graph $G(n,p)$ a.s. contains a Hamiltonian path. We believe that a more general 
result should be true, i.e.,  such a random graph should 
already contain a.s. every 
tree on $n$ vertices with maximum degree at most $d$. Our methods are insufficient to attack this problem.
Using Theorem \ref{random} we can only prove the following much weaker result.
Let $T$ be a tree on $n$ vertices with at least $\epsilon n$ leaves, then there exists a constant
$a(\epsilon,d)$ such that the random graph 
$G\big(n,\frac{a\log n}{n}\big)$ a.s. contains $T$. Here is a
brief sketch of the proof: we split $G(n,p)$ into two random 
graphs $G(n,p_1)$ and $G(n,p_2)$, where $1-p=(1-p_1)(1-p_2)$, 
$p_1=\Theta(1/n)$, $p_2=\Theta(\log n/n)$. Let $T'$ be the tree
obtained from $T$ by deleting its leaves. We use Theorem
\ref{random} to embed a copy of $T'$ in $G(n,p_1)$. Then we
expose the edges of $G(n,p_2)$ between the set of vertices $V_0$
of $G$, not occupied by a copy of $T'$, and the rest of the 
graph, and embed the leaves of $T$ in $V_0$ using matching-type
results. 

\item
Besides the model $G(n,p)$, another model of random
graphs, drawing a lot of attention is the model of random
regular graphs. A random regular graph $G_{n,D}$ is obtained by sampling uniformly
at random over the set of all simple $D$-regular graphs on a fixed set of $n$ vertices. 
By the result of Friedman, Kahn and Szemer\'edi \cite{FKS} the second eigenvalue of $G_{n,D}$ is almost
surely at most $O(\sqrt{D})$ (see \cite{Fr} for a more precise result). Therefore our 
Theorem \ref{pseudo-random} immediately implies that if 
$D=D(d,\epsilon)$ is sufficiently large then $G_{n,D}$ a.s. contains every tree on
$(1-\epsilon)n$ vertices with maximum degree $d$. 

This result as well as the result of Theorem \ref{pseudo-random} are probably not
optimal. We suspect that sufficiently large spectral gap (as a 
function of $d$ only) already suffices 
to guarantee the embedding of every spanning tree of 
bounded degree in a graph $G$ of order $n$.
This is not known even for the Hamiltonian path, and the best 
result in this case, obtained in \cite{KS},
requires spectral gap of order  roughly $\log n$.
\end{itemize}

\vspace{0.15cm}
\noindent
{\bf Acknowledgment.}\, A major part of this work was carried
out when the authors were visiting Microsoft Research at Redmond, WA.
We would like to thank the members of the Theory Group at
Microsoft Research for their hospitality and for creating a
stimulating research environment.

\end{document}